\documentclass[11pt]{article}

\usepackage{latexsym}
\usepackage{amssymb}

\newcommand{\Z}{{\mathbb Z}} %
\newcommand{\R}{{\mathbb R}} %
\newcommand{\CC}{{\mathbb C}} %
\newcommand{\N}{{\mathbb N}} %
\newcommand{\PP}{{\mathbb P}}%

\newtheorem{thm}{Theorem}[section]

\newcommand{\ld}{\lambda}

\newcommand{\img}{\mathop{\rm Im}\nolimits}

\newcommand{\Ker}{\mathop{\rm Ker}\nolimits}

\newcommand{\Lin}{\mathop{\rm Lin}\nolimits}
\newcommand{\mspan}{\mathop{\rm span}\nolimits}

\begin{document}

\begin{center}
{\Large\bf
Truncated matricial moment problems on a finite interval: the operator approach.}
\end{center}
\begin{center}
{\bf S.M. Zagorodnyuk}
\end{center}

\section{Introduction.}
In this paper we study the following problem:
to find a non-decreasing matrix function
$M(x) = ( m_{k,l}(x) )_{k,l=0}^{N-1}$ on $[a,b]$, which is left-continuous in $(a,b)$, $M(a)=0$, such that
\begin{equation}
\label{f1_1}
\int_a^b x^n dM(x) = S_n,\qquad n=0,1,...,l,
\end{equation}
where $\{ S_n \}_{n=0}^l$ is a given sequence of Hermitian $(N\times N)$ complex matrices, $N\in\N$, $l\in\Z_+$.
Here $a,b\in\R$: $a<b$.

\noindent
In the scalar case this problem was solved by M.G.~Krein, see~\cite{cit_1000_KN}.
Recently, a deep investigation of the matrix moment problem~(\ref{f1_1}) was completed by
A.E.~Choque Rivero, Yu.M.~Dyukarev, B.~Fritzsche and B.~Kirstein, see~\cite{cit_2000_CDFK},\cite{cit_3000_CDFK}.
These authors used the Potapov method for interpolating problems which was enriched by the Sachnovich method of
operator identities.

\noindent
Set
\begin{equation}
\label{f1_2}
\Gamma_k = ( S_{i+j} )_{i,j=0}^k = \left(
\begin{array}{cccc} S_0 & S_1 & \ldots & S_k\\
S_1 & S_2 & \ldots & S_{k+1}\\
\vdots & \vdots & \ddots & \vdots\\
S_k & S_{k+1} & \ldots & S_{2k}\end{array}
\right),\qquad k\in\Z_+:\ 2k\leq l;
\end{equation}
\begin{equation}
\label{f1_3}
\widetilde\Gamma_k = ( -ab S_{i+j} + (a+b) S_{i+j+1} - S_{i+j+2} )_{i,j=0}^{k-1},\qquad k\in\Z_+:\ 2k\leq l.
\end{equation}
If we choose an arbitrary element $f =(f_0,f_1,\ldots,f_{N-1})$, where all $f_k$ are some polynomials
and calculate
$\int_a^b f dM f^*$, one can easily deduce that
\begin{equation}
\label{f1_4}
\Gamma_k \geq 0,\qquad k\in\Z_+:\ 2k\leq l.
\end{equation}
In the case of an odd number of prescribed moments $l=2d$, the strong result of
A.E.~Choque Rivero, Yu.M.~Dyukarev, B.~Fritzsche and B.~Kirstein is that conditions
\begin{equation}
\label{f1_5}
\Gamma_d \geq 0,\  \widetilde\Gamma_d \geq 0,
\end{equation}
are necessary and sufficient for the solvability of the matrix moment problem~(\ref{f1_1}),
see~\cite[Theorem 1.3, p. 106]{cit_3000_CDFK}.
For the case $\Gamma_d>0,\widetilde\Gamma_d>0$, they parameterized all solutions of the moment
problem via a linear fractional transformation where the set of parameters consisted of some
distinguished pairs of meromorphic matrix-valued functions.

\noindent
Set
\begin{equation}
\label{f1_6}
H_k = ( -a S_{i+j} + S_{i+j+1} )_{i,j=0}^k,\
\widetilde H_k = ( b S_{i+j} - S_{i+j+1} )_{i,j=0}^k,\quad  k\in\Z_+:\ 2k+1\leq l.
\end{equation}
In the case $l=2d+1$, the analogous to the result for $l=2d$, the result of
A.E.~Choque Rivero, Yu.M.~Dyukarev, B.~Fritzsche and B.~Kirstein states that conditions
\begin{equation}
\label{f1_7}
H_d \geq 0,\  \widetilde H_d \geq 0,
\end{equation}
are necessary and sufficient for the solvability of the matrix moment problem~(\ref{f1_1}),
see~\cite[Theorem 1.3, p. 127]{cit_2000_CDFK}.
For the case $H_d>0,\widetilde H_d>0$, they parameterized all solutions of the moment
problem via a linear fractional transformation. The set of parameters consisted of some
distinguished pairs of meromorphic matrix-valued functions.

\noindent
In this work we will study the matrix moment problem~(\ref{f1_1}) by virtue of the operator approach based on
the use of the generalized resolvents of some symmetric operators. In the study of the
classical Hamburger moment problem this approach finds its origin in the papers of
M.A.~Neumark~\cite{cit_4000_N},\cite{cit_5000_N} and
M.G.~Krein, M.A. Krasnoselskiy~\cite{cit_6000_KK}, see also~\cite{cit_7000_A}.
All these authors used orthogonal polynomials connected with a Jacobi matrix related to the moment problem.
Lately, we showed that this method in a general setting can be applied to the Hamburger moment problem
both in the non-degenerate and degenerate cases, see~\cite{cit_8000_Z}.
Our goal here is to describe all solutions of the matrix moment problem~(\ref{f1_1}) in a general case. This
means that no conditions besides solvability of the moment problem will be assumed.
At first, we study the case of an odd number of prescribed moments $l=2d$, and then we shall reduce
the case of an even number of moments $l=2d+1$ to the previous case ($d\in\Z_+$).
In our study we shall use the basic results of M.G.~Krein and I.E.~Ovcharenko on generalized sc-resolvents
of symmetric contractions, as well as M.G.~Krein's theory of self-adjoint extensions of semi-bounded
symmetric operators, see~\cite{cit_9000_KO},\cite{cit_10000_KO},\cite{cit_11000_K}.

\noindent
{\bf Notations.}  As usual, we denote by $\R, \CC, \N, \Z, \Z_+$
the sets of real, complex, positive integer, integer, non-negative integer numbers,
respectively. The space of $n$-dimensional complex vectors $a = (a_0,a_1,\ldots,a_{n-1})$, will be
denoted by $\CC^n$, $n\in\N$.
If $a\in \CC^n$ then $a^*$ means the complex conjugate vector.
For a complex $(n\times n)$ matrix $A$, we denote by $\Ker A$ a set $\{ x\in\CC^n:\ Ax=0 \}$.
By $\PP$ we denote a set of all complex polynomials and by $\PP_d$ we mean all complex polynomials with
degrees less or equal to $d$, $d\in\Z_+$, (including the zero polynomial).
Let $M(x)$ be a left-continuous non-decreasing matrix function $M(x) = ( m_{k,l}(x) )_{k,l=0}^{N-1}$
on $\R$, $M(-\infty)=0$, and $\tau_M (x) := \sum_{k=0}^{N-1} m_{k,k} (x)$;
$\Psi(x) = ( dm_{k,l}/ d\tau_M )_{k,l=0}^{N-1}$.  We denote by $L^2(M)$ a set (of classes of equivalence)
of vector functions $f: \R\rightarrow \CC^N$, $f = (f_0,f_1,\ldots,f_{N-1})$, such
that (see, e.g.,~\cite{cit_11500_MM})
$$ \| f \|^2_{L^2(M)} := \int_\R  f(x) \Psi(x) f^*(x) d\tau_M (x) < \infty. $$
The space $L^2(M)$ is a Hilbert space with the scalar product
$$ ( f,g )_{L^2(M)} := \int_\R  f(x) \Psi(x) g^*(x) d\tau_M (x),\qquad f,g\in L^2(M). $$
For a separable Hilbert space $H$ we denote by $(\cdot,\cdot)_H$ and $\| \cdot \|_H$ the scalar
product and the norm in $H$, respectively. The indices may be omitted in obvious cases.

\noindent
For a linear operator $A$ in $H$ we denote by $D(A)$ its domain, by $R(A)$ its range, and by
$A^*$ we denote its adjoint if it exists. If $A$ is bounded, then $\| A \|$ stands for its operator norm.
For a set of elements $\{ x_n \}_{n\in B}$ in $H$, we
denote by $\Lin\{ x_n \}_{n\in B}$ and $\mspan\{ x_n \}_{n\in B}$ the linear span and the closed
linear span (in the norm of $H$), respectively, where $B$ is an arbitrary set of indices.
For a set $M\subseteq H$ we denote by $\overline{M}$ the closure of $M$ with respect to the norm of $H$.
By $E_H$ we denote the identity operator in $H$, i.e. $E_H x = x$, $x\in H$.
If $H_1$ is a subspace of $H$, by $P_{H_1} = P_{H_1}^{H}$ we denote the operator of the orthogonal projection on $H_1$
in $H$. A set of linear bounded operators which map $H$ into $H$ we denote by $[H]$.

\section{The case of an odd number of given moments: solvability and a description  of solutions.}
We shall use the following important fact (see, e.g., \cite[p.215]{cit_12000}):
\begin{thm}
\label{t2_1}
Let $K = (K_{n,m})_{n,m=0}^{r} \geq 0$ be a positive semi-definite complex $((r+1)\times(r+1))$ matrix, $r\in\Z_+$.
Then there exist a finite-dimensional Hilbert space $H$ with a scalar product $(\cdot,\cdot)$ and
a sequence $\{ x_n \}_{n=0}^r$ in $H$, such that
\begin{equation}
\label{f2_1}
K_{n,m} = (x_n,x_m),\qquad n,m=0,1,...,r,
\end{equation}
and $\mspan\{ x_n \}_{n=0}^r = H$.
\end{thm}
P r o o f.
Let $\{ x_n \}_{n=0}^r$ be an arbitrary orthonormal basis in $\CC^n$.
Introduce the following functional:
\begin{equation}
\label{f2_2}
[x,y] = \sum_{n,m=0}^r K_{n,m} a_n\overline{b_m},
\end{equation}
for $x,y\in \CC^n$,
$$ x=\sum_{n=0}^r a_n x_n,\quad y=\sum_{m=0}^r b_m x_m,\quad a_n,b_m\in\CC. $$
The space $\CC^n$ equipped with $[\cdot,\cdot]$ will be a quasi-Hilbert space.
Factorizing and making the completion we obtain the  required space $H$
(see, e.g.,~\cite[p. 10-11]{cit_13000}).
$\Box$

\noindent
Consider the matrix moment problem~(\ref{f1_1}) with $l=2d$, $d\in\N$.
Suppose that $\Gamma_d\geq 0$ (as we noticed in the Introduction,
condition~(\ref{f1_4}) is necessary for the solvability of the moment problem).
Let
$\Gamma_d = (\gamma_{d;n,m})_{n,m=0}^{(d+1)N-1}$, $\gamma_{d;n,m}\in\CC$.
By Theorem~\ref{t2_1} there exist a finite-dimensional Hilbert space $H$ and
a sequence $\{ x_n \}_{n=0}^{ (d+1)N-1 }$ in $H$, such that
\begin{equation}
\label{f2_3}
(x_n,x_m) = \gamma_{d;n,m},\qquad n,m=0,1,...,(d+1)N-1,
\end{equation}
and
$\mspan\{ x_n \}_{n=0}^{ (d+1)N-1 } = \Lin\{ x_n \}_{n=0}^{ (d+1)N-1 } = H$.
Notice that
\begin{equation}
\label{f2_3_1}
\gamma_{d;rN+j,tN+n} = s_{r+t}^{j,n},\qquad 0\leq j,n \leq N-1;\quad 0\leq r,t\leq d,
\end{equation}
where
$$ S_n = ( s_n^{k,l} )_{k,l=0}^{N-1},\qquad n\in\Z_+, $$
are the given moments.
From~(\ref{f2_3_1}) it follows that
\begin{equation}
\label{f2_3_2}
\gamma_{d;a+N,b} = \gamma_{d;a,b+N},\qquad a=rN+j,\ b=tN+n,\ 0\leq j,n \leq N-1;\quad 0\leq r,t\leq d-1.
\end{equation}
In fact, we can write
$$ \gamma_{d;a+N,b} = \gamma_{d;(r+1)N+j,tN+n} = s_{r+t+1}^{j,n} = \gamma_{d;rN+j,(t+1)N+n} =
\gamma_{d;a,b+N}. $$
Set $H_a = \{ x_n \}_{n=0}^{ dN-1 }$. We introduce the following operator:
\begin{equation}
\label{f2_4}
A x = \sum_{k=0}^{dN-1} \alpha_k x_{k+N},\qquad x\in H_a,\ x=\sum_{k=0}^{dN-1}\alpha_k x_k.
\end{equation}
The following proposition shows when the operator $A$ is correctly defined.
\begin{thm}
\label{t2_2}
Let a matrix moment problem~(\ref{f1_1}) with $l=2d$, $d\in\N$, be given.
The moment problem has a solution if and only if
conditions~(\ref{f1_5}) are true and
\begin{equation}
\label{f2_5}
\Ker \Gamma_{d-1} \subseteq \Ker \widehat\Gamma_{d-1},
\end{equation}
where $\widehat\Gamma_{d-1} = (S_{i+j+2})_{i,j=0}^{d-1}$.

If conditions~(\ref{f1_5}),(\ref{f2_5}) are satisfied then the operator $A$ in~(\ref{f2_4})
is correctly defined and the following operator:
\begin{equation}
\label{f2_6}
Bx = \frac{2}{b-a} A - \frac{a+b}{b-a} E_H,\qquad x\in H_a,
\end{equation}
is a contraction in $H$ (i.e. $\| B \| \leq 1$). Moreover, operators $A$ and $B$ are Hermitian.
\end{thm}
P r o o f.
Let the matrix moment problem~(\ref{f1_1}) has a solution
$M(x) = (m_{k,l}(x) )_{k,l=0}^{N-1}$.
Consider the space $L^2(M)$ and
let $Q$ be the operator of multiplication by an independent variable in $L^2(M)$.
The operator $Q$ is self-adjoint and its resolution of unity is (see~\cite{cit_11500_MM})
\begin{equation}
\label{f2_7}
E_b - E_a = E([a,b)): h(x) \rightarrow \chi_{[a,b)}(x) h(x),
\end{equation}
where $\chi_{[a,b)}(x)$ is the characteristic function of an interval $[a,b)$, $-\infty\leq a<b\leq +\infty$.

\noindent
Set $\vec e_k = (e_{k,0},e_{k,1},\ldots,e_{k,N-1})$, $e_{k,j}=\delta_{k,j}$, $0\leq j\leq N-1$,
for $k=0,1,\ldots N-1$.
A set of (classes of equivalence of) functions $f\in L^2(M)$ such that
(the corresponding class includes) $f=(f_0,f_1,\ldots, f_{N-1})$, $f_j\in\PP_d$, we denote
by $\PP^2_d(M)$ and call a set of vector polynomials of order $d$ in $L^2(M)$.
Set $L^2_{d,0}(M) = \overline{ \PP^2_d(M) }$. Since $\PP^2_d(M)$ is finite-dimensional, we have
$L^2_{d,0}(M) = \PP^2_d(M)$.

For an arbitrary polynomial (in a class) from $\PP^2_d(M)$ there exists a
unique representation of the following form:
\begin{equation}
\label{f2_8}
f(x) = \sum_{k=0}^{N-1} \sum_{j=0}^d \alpha_{k,j} x^j \vec e_k,\quad \alpha_{k,j}\in\CC.
\end{equation}
Let a polynomial $g\in \PP^2_d(M)$ have a representation
\begin{equation}
\label{f2_9}
g(x) = \sum_{l=0}^{N-1} \sum_{r=0}^d \beta_{l,r} x^r \vec e_l,\quad \beta_{l,r}\in\CC.
\end{equation}
We can write
$$ (f,g)_{L^2(M)} = \sum_{k,l=0}^{N-1} \sum_{j,r=0}^d \alpha_{k,j}\overline{\beta_{l,r}}
\int_\R x^{j+r} \vec e_k dM(x) \vec e_l^* = \sum_{k,l=0}^{N-1}
\sum_{j,r=0}^d \alpha_{k,j}\overline{\beta_{l,r}} * $$
\begin{equation}
\label{f2_10}
* \int_\R x^{j+r} dm_{k,l}(x) = \sum_{k,l=0}^{N-1} \sum_{j,r=0}^d \alpha_{k,j}\overline{\beta_{l,r}}
s_{j+r}^{k,l},
\end{equation}
where
\begin{equation}
\label{f2_11}
S_n = ( s_n^{k,l} )_{k,l=0}^{N-1},\qquad n\in\Z_+,
\end{equation}
are the given moments.
On the other hand, we can write
$$ \left( \sum_{j=0}^d \sum_{k=0}^{N-1} \alpha_{k,j} x_{jN+k},
\sum_{r=0}^d \sum_{l=0}^{N-1} \beta_{l,r} x_{rN+l} \right)_H =
\sum_{k,l=0}^{N-1} \sum_{j,r=0}^d \alpha_{k,j}\overline{\beta_{l,r}}
(x_{jN+k}, x_{rN+l})_H  = $$
\begin{equation}
\label{f2_12}
= \sum_{k,l=0}^{N-1} \sum_{j,r=0}^d \alpha_{k,j}\overline{\beta_{l,r}}
\gamma_{d;jN+k,rN+l}
= \sum_{k,l=0}^{N-1} \sum_{j,r=0}^d \alpha_{k,j}\overline{\beta_{l,r}}
s_{j+r}^{k,l}.
\end{equation}
From relations~(\ref{f2_10}),(\ref{f2_12}) it follows that
\begin{equation}
\label{f2_13}
(f,g)_{L^2(M)} = \left( \sum_{j=0}^d \sum_{k=0}^{N-1} \alpha_{k,j} x_{jN+k},
\sum_{r=0}^d \sum_{l=0}^{N-1} \beta_{l,r} x_{rN+l} \right)_H.
\end{equation}
Set
\begin{equation}
\label{f2_14}
Vf = \sum_{j=0}^d \sum_{k=0}^{N-1} \alpha_{k,j} x_{jN+k},
\end{equation}
for $f(x)\in \PP^2_d(M)$, $f(x) = \sum_{k=0}^{N-1} \sum_{j=0}^d \alpha_{k,j} x^j \vec e_k$,
$\alpha_{k,j}\in\CC$.

\noindent
If $f$, $g$ have representations~(\ref{f2_8}),(\ref{f2_9}), and $\| f-g \|_{L^2(M)} = 0$, then
from~(\ref{f2_13}) it follows that
$$ \| Vf - Vg \|_H^2 = (V(f-g),V(f-g))_H = ( f-g,f-g )_{L^2(M)} = \| f-g\|_{L^2(M)}^2 = 0.$$
Thus, $V$ is a correctly defined operator from $\PP^2_d(M)$ to $H$.

\noindent
Relation~(\ref{f2_13}) shows that $V$ is an isometric transformation from $\PP^2_d(M)$ onto
$\Lin\{ x_n \}_{n=0}^{d(N+1)-1}$.
Thus, $V$ is an isometric transformation from $L^2_{d,0}(M)$ onto $H$.
In particular, we note that
\begin{equation}
\label{f2_15}
V x^j \vec e_k = x_{jN+k},\qquad 0\leq j\leq d;\quad 0\leq k\leq N-1.
\end{equation}
Set $L^2_{d,1} (M) := L^2(M)\ominus L^2_{d,0} (M)$, and
$U := V\oplus E_{L^2_{d,1} (M)}$. The operator $U$ is
an isometric transformation from $L^2(M)$ onto $H\oplus L^2_{d,1} (M)=:\widehat H$.
Set
\begin{equation}
\label{f2_16}
\widehat A := UQU^{-1}.
\end{equation}
The operator $\widehat A$ is a self-adjoint operator in $\widehat H$.
Notice that
$$ UQU^{-1} x_{jN+k} = VQV^{-1} x_{jN+k} = VQ x^j \vec e_k = V x^{j+1} \vec e_k =
x_{(j+1)N+k} = x_{jN+k+N},$$
$$ 0\leq j\leq d-1;\quad 0\leq k\leq N-1. $$
By linearity we get
$$ UQU^{-1} x = \sum_{k=0}^{dN-1} \alpha_k x_{k+N},\qquad x\in H_a,\ x=\sum_{k=0}^{dN-1}\alpha_k x_k. $$
Consequently, the operator $A$ in~(\ref{f2_4}) is correctly defined and
\begin{equation}
\label{f2_17}
A = \widehat A|_{H_a}.
\end{equation}
Since $A$ is correctly defined, from the equality
\begin{equation}
\label{f2_18}
\sum_{k=0}^{dN-1} \xi_k x_k = 0,
\end{equation}
with some complex numbers $\xi_k$, it should follow the equality
\begin{equation}
\label{f2_19}
\sum_{k=0}^{dN-1} \xi_k x_{k+N} = 0.
\end{equation}
On the other hand, the equality~(\ref{f2_18}) is equivalent to the equalities
\begin{equation}
\label{f2_20}
\sum_{k=0}^{dN-1} \xi_k (x_k,x_l) =  \sum_{k=0}^{dN-1} \xi_k \gamma_{d;k,l} = 0,\quad l=0,1,...,dN-1.
\end{equation}
Analogously, the equality~(\ref{f2_19}) is equivalent to the equalities
\begin{equation}
\label{f2_21}
\sum_{k=0}^{dN-1} \xi_k (x_{k+N},x_{l+N}) =  \sum_{k=0}^{dN-1} \xi_k \gamma_{d;k+N,l+N} = 0,\quad l=0,1,...,dN-1.
\end{equation}
If we shall use the matrix notations, the equality
\begin{equation}
\label{f2_22}
(\xi_0,\xi_1,...,\xi_{dN-1}) (\gamma_{d;k,l})_{k,l=0}^{dN-1} = 0,
\end{equation}
implies the equality
\begin{equation}
\label{f2_23}
(\xi_0,\xi_1,...,\xi_{dN-1}) (\gamma_{d;k+N,l+N})_{k,l=0}^{dN-1} = 0.
\end{equation}
Thus, relation~(\ref{f2_5}) is true.

\noindent
Consider the following operators:
\begin{equation}
\label{f2_24}
R := \frac{2}{b-a} Q - \frac{a+b}{b-a} E_{L^2(M)},
\end{equation}
\begin{equation}
\label{f2_25}
\widehat B := URU^{-1} = \frac{2}{b-a} \widehat A - \frac{a+b}{b-a} E_{\widehat H}.
\end{equation}
Define an operator $B$ by the equality~(\ref{f2_6}). From~(\ref{f2_17}),(\ref{f2_25}) we get
\begin{equation}
\label{f2_26}
B = \widehat B|_{H_a}.
\end{equation}
For an arbitrary $f\in D(R)=D(Q)$ we can write
$$ \| Rf \|^2_{L^2(M)} = \int_a^b \left| \frac{2}{b-a} x - \frac{a+b}{b-a} \right|^2
f(x) dM(x) f^*(x) \leq \int_a^b f(x) dM(x) f^*(x) = $$
$$ = \| f \|^2, $$
and therefore the operators $R$, $\widehat B$ and $B$ are contractions.
Since $R$ is Hermitian, the operators $\widehat B$, $B$ are Hermitian, as well.
Choose an arbitrary $x\in H_a$, $x=\sum_{k=0}^{dN-1} \alpha_k x_k$, and write
$$ 0 \leq \| x \|^2 - \| Bx \|^2 = \sum_{k,j=0}^{dN-1} \alpha_k \overline{\alpha_j} (x_k,x_j) -
\sum_{k,j=0}^{dN-1} \alpha_k \overline{\alpha_j} (Bx_k,Bx_j) = $$
$$ = \sum_{k,j=0}^{dN-1} \alpha_k \overline{\alpha_j} \gamma_{d;k,j} -
\sum_{k,j=0}^{dN-1} \alpha_k \overline{\alpha_j} \left( \frac{2}{b-a} x_{k+N} - \frac{a+b}{b-a} x_k,
\frac{2}{b-a} x_{j+N} - \frac{a+b}{b-a} x_j \right) = $$
$$ = \sum_{k,j=0}^{dN-1} \alpha_k \overline{\alpha_j}
\left( \gamma_{d;k,j} - \frac{4}{(b-a)^2} \gamma_{d;k+N,j+N} + \frac{2(a+b)}{(b-a)^2} \gamma_{d;k+N,j} +
\frac{2(a+b)}{(b-a)^2} \gamma_{d;k,j+N} - \right. $$
\begin{equation}
\label{f2_27}
\left. - \frac{(a+b)^2}{(b-a)^2} \gamma_{d;k,j} \right).
\end{equation}
If we multiply the both sides of the latter inequality by $(b-a)^2$ and use~(\ref{f2_3_2}) we
get
$$ 0 \leq \sum_{k,j=0}^{dN-1} \alpha_k \overline{\alpha_j}
\left( (b-a)^2 \gamma_{d;k,j} - 4\gamma_{d;k+N,j+N} + 4(a+b) \gamma_{d;k+N,j} -
(a+b)^2 \gamma_{d;k,j} \right) = $$
$$ = \sum_{k,j=0}^{dN-1} \alpha_k \overline{\alpha_j}
( -4ab \gamma_{d;k,j} - 4\gamma_{d;k+N,j+N} + 4(a+b) \gamma_{d;k+N,j}).  $$
Therefore $0\leq (-ab \gamma_{d;k,j} - \gamma_{d;k+N,j+N} + (a+b) \gamma_{d;k+N,j})_{k,j=0}^{dN-1} =
\widetilde\Gamma_{d}$, and the second relation in~(\ref{f1_5}) is true.
The necessity of the first relation in~(\ref{f1_5}) was discussed in the Introduction.
Thus, we established the necessity of conditions~(\ref{f1_5}),(\ref{f2_5}) for the solvability of the
moment problem~(\ref{f1_1}).

Let a matrix moment problem~(\ref{f1_1}) be given and conditions~(\ref{f1_5}),(\ref{f2_5}) be true.
We construct a Hilbert space $H$ and a sequence $\{ x_n \}_{n=0}^{ (d+1)N-1 }$ in $H$, such that
relation~(\ref{f2_3}) is true and $\mspan\{ x_n \}_{n=0}^{ (d+1)N-1 } = H$.
Condition~(\ref{f2_5}) provides that equality~(\ref{f2_18}) will imply equality~(\ref{f2_19}).
This means that the operator $A$ in~(\ref{f2_4}) is correctly defined.
We define an operator $B$ by the equality~(\ref{f2_6}).
For arbitrary $x,y\in H_a = \mspan\{ x_n \}_{n=0}^{ dN-1 }$, $x = \sum_{k=0}^{dN-1} \alpha_k x_k$,
$y = \sum_{j=0}^{dN-1} \beta_j x_j$, using~(\ref{f2_3_2}) we can write
$$ (Ax,y) = \sum_{k,j=0}^{dN-1} \alpha_k \overline{\beta_j} (x_{k+N},x_j) =
\sum_{k,j=0}^{dN-1} \alpha_k \overline{\beta_j} (x_{k},x_{j+N}) = (x,Ay). $$
Thus, operators $A$ and $B$ are Hermitian.
If we use relation~(\ref{f2_27}) (except the first inequality in it) and the second condition
in~(\ref{f1_5}), we obtain that $B$ is a contraction.
By Krein's theorem~\cite[Theorem 2, p. 440]{cit_11000_K}, there exists
a self-adjoint extension $\widetilde B$ of the operator $B$ in $H$ with the same norm as $B$
(and therefore it is a contraction).
Let
\begin{equation}
\label{f2_28}
\widetilde B = \int_{-1}^1 \ld d\widetilde E_\ld,
\end{equation}
where $\{ \widetilde E_\ld \}$ be the left-continuous in $[-1,1)$, right-continuous
at the point $1$, constant outside $[-1,1]$, orthogonal resolution of unity of
$\widetilde B$.
Choose an arbitrary $\alpha$, $0\leq \alpha\leq d(N+1)-1$, $\alpha=rN + j$, $0\leq r\leq d$, $0\leq j\leq N-1$.
Notice that
$$ x_\alpha = x_{rN+j} = A x_{(r-1)N+j} = ... = A^r x_j. $$
Then choose an arbitrary $\beta$, $0\leq \beta\leq d(N+1)-1$, $\beta=tN + n$, $0\leq t\leq d$, $0\leq n\leq N-1$.
Using~(\ref{f2_3_1}) we can write
$$ s_{r+t}^{j,n} = \gamma_{d;rN+j,tN+n} = ( x_{rN+j},x_{tN+n} )_H = (A^r x_j, A^t x_n)_H = $$
$$ =
\left( \left( \frac{b-a}{2} B + \frac{a+b}{2} E_H \right)^r x_j,
  \left( \frac{b-a}{2} B + \frac{a+b}{2} E_H \right)^t x_n \right)_H
= $$
$$ =
\left( \left( \frac{b-a}{2} \widetilde B + \frac{a+b}{2} E_H \right)^{r+t} x_j, x_n \right)_H
= \int_{-1}^1 \left( \frac{b-a}{2} \ld + \frac{a+b}{2} \right)^{r+t} d( \widetilde E_\ld x_j,x_n )_H. $$
Set
$$ \widetilde m_{j,n}(x) = ( \widetilde E_{ \frac{2}{b-a}x - \frac{a+b}{b-a} } x_j,x_n )_H,\qquad
0\leq j,n\leq N-1. $$
Then
\begin{equation}
\label{f2_29}
s_{r+t}^{j,n} = \int_{a}^b x^{r+t}
d \widetilde m_{j,n}(x),\qquad 0\leq j,n\leq N-1,\ 0\leq r,t\leq d.
\end{equation}
From~relation~(\ref{f2_29}) we derive that the matrix function
$\widetilde M(\ld)=(\widetilde m_{j,n}(x))_{j,n=0}^{N-1}$ is
a solution of the matrix Hamburger moment  problem~(\ref{f1_1}) (Properties of the
orthogonal resolution of unity provide that $\widetilde M (\ld)$ is left-continuous in $(a,b)$, non-decreasing and
$\widetilde M(a) = 0$).

It remains to prove the last statement of the Theorem.
If conditions~(\ref{f1_5}),(\ref{f2_5}) are satisfied then we proved that the moment
problem~(\ref{f1_1}) has a solution. In this case we showed that the operator $A$ in~(\ref{f2_4})
is correctly defined and the operator $B$ in~(\ref{f2_6})
is a Hermitian contraction in $H$. The fact that operators $A$ and $B$ are Hermitian was established, as well.
$\Box$

We shall continue our considerations before the statement of Theorem~\ref{t2_2}. We assume that
conditions~(\ref{f1_5}),(\ref{f2_5}) are true. Therefore the operators $A$ in~(\ref{f2_4}) and
$B$ in~(\ref{f2_6}) are correctly defined Hermitian operators and $\| B \|\leq 1$.

\noindent
Let $\widehat B$ be an arbitrary self-adjoint extension of $B$ in a Hilbert space $\widehat H\supseteq H$.
Let $R_z(\widehat B)$ be the resolvent of $\widehat B$ and $\{ \widehat E_\lambda\}_{\ld\in\R}$
be an orthogonal resolution of unity of $\widehat B$. Recall that the operator-valued function
$\mathbf R_z = P_H^{\widehat H} R_z(\widehat B)$ is called a generalized resolvent of $B$, $z\in\CC\backslash\R$.
The function
$\mathbf E_\lambda = P_H^{\widehat H} \widehat E_\ld$, $\ld\in\R$, is a spectral
function of a symmetric operator $B$ (e.g.~\cite{cit_13500_S}).
There exists a one-to-one correspondence between generalized resolvents and (left-continuous or
normalized in another way) spectral functions
established by the following relation (\cite{cit_14000_AG}):
\begin{equation}
\label{f2_30}
(\mathbf R_z f,g)_H = \int_\R \frac{1}{\ld - z} d( \mathbf E_\ld f,g)_H,\qquad f,g\in H,\ z\in\CC\backslash\R.
\end{equation}
To obtain the spectral function from relation~(\ref{f2_30}), one should use the Stieltjes-Perron
inversion formula (e.g.~\cite{cit_7000_A}).

\noindent
In the case when $\widehat B$ is a self-adjoint contraction,
the generalized resolvent $\mathbf R_z = P_H^{\widehat H} R_z(\widehat B)$ is called a generalized sc-resolvent of $B$,
see~\cite{cit_9000_KO},\cite{cit_10000_KO}.
The corresponding spectral function of $B$ we shall call a sc-spectral function of $B$.

\noindent
Let
\begin{equation}
\label{f2_31}
\widehat B = \int_{-1}^1 \ld d\widehat E_\ld,
\end{equation}
where $\{ \widehat E_\ld \}$ be the left-continuous in $[-1,1)$, right-continuous
at the point $1$, constant outside $[-1,1]$, orthogonal resolution of unity of
$\widehat B$.
In a similar manner as after~(\ref{f2_28}) we obtain that
\begin{equation}
\label{f2_32}
s_{r+t}^{j,n} = \int_{a}^b x^{r+t}
d \widehat m_{j,n}(x),\qquad 0\leq j,n\leq N-1,\ 0\leq r,t\leq d,
\end{equation}
where
\begin{equation}
\label{f2_33}
\widehat m_{j,n}(x) = ( P^{\widehat H}_H \widehat E_{ \frac{2}{b-a}x - \frac{a+b}{b-a} } x_j,x_n )_H,\qquad
0\leq j,n\leq N-1.
\end{equation}
Thus, the function $\widehat M(x) = (\widehat m_{j,n}(x))_{j,n=0}^{N-1}$ is a solution of the
matrix moment problem~(\ref{f1_1}).
\begin{thm}
\label{t2_3}
Let a matrix moment problem~(\ref{f1_1}) with $l=2d$, $d\in\N$, be given.
Suppose that conditions~(\ref{f1_5}),(\ref{f2_5}) are true.
All solutions of the moment problem have the following form
\begin{equation}
\label{f2_34}
M(x) = (m_{j,n}(x))_{j,n=0}^{N-1},\quad
m_{j,n}(x) = ( \mathbf E_{ \frac{2}{b-a}x - \frac{a+b}{b-a} } x_j,x_n )_H,\qquad
0\leq j,n\leq N-1,
\end{equation}
where $\mathbf E_z$ is a left-continuous in $[-1,1)$, right-continuous
at the point $1$, constant outside $[-1,1]$ sc-spectral function of the operator $B$ defined in~(\ref{f2_6}).

\noindent
Moreover, the correspondence between all solutions of the moment problem and
left-continuous in $[-1,1)$, right-continuous at the point $1$, constant outside $[-1,1]$
sc-spectral functions
of $B$ in~(\ref{f2_34}) is one-to-one.
\end{thm}
P r o o f.
Choose an arbitrary left-continuous in $[-1,1)$, right-continuous at the point $1$, constant outside $[-1,1]$
sc-spectral function $\mathbf E_z$ of the operator $B$ from~(\ref{f2_6}).
This function corresponds to a left-continuous in $[-1,1)$, right-continuous at the point $1$, constant
outside $[-1,1]$
resolution of unity $\{ \widehat E_\ld \}$ of a self-adjoint contraction $\widehat B\supseteq B$
in a Hilbert space $\widehat H\supseteq H$. Considerations before the statement of the Theorem show
that formula~(\ref{f2_34}) defines a solution of the moment problem.

\noindent
On the other hand, let $M(x)=(m_{k,l}(x) )_{k,l=0}^{N-1}$ be an arbitrary solution of the matrix
moment problem~(\ref{f1_1}). Proceeding like at the beginning of the Proof of Theorem~\ref{t2_2},
we shall construct a self-adjoint contraction $\widehat B\supseteq B$ in a  space $\widehat H\supseteq H$.
Repeating arguments before the statement of the last theorem, we obtain that
the function $\widehat M(x) = (\widehat m_{j,n}(x))_{j,n=0}^{N-1}$, where
$\widehat m_{j,n}(x)$ are given by~(\ref{f2_33}), is a solution of the moment problem.

\noindent
Choose an arbitrary $z\in\CC\backslash\R$ and write
$$ \int_{-1}^1 \frac{1}{\ld - z} d( \widehat E_\ld x_k, x_j)_{\widehat H} =
   \left( \int_{-1}^1 \frac{1}{\ld - z} d\widehat E_\ld x_k, x_j \right)_{\widehat H} = $$
$$ = \left( U^{-1} \int_{-1}^1 \frac{1}{\ld - z} d\widehat E_\ld x_k, U^{-1} x_j \right)_{L^2(M)}
   = \left( \int_{-1}^1 \frac{1}{\ld - z} d U^{-1} \widehat E_\ld U \vec e_k, \vec e_j \right)_{L^2(M)} = $$
$$ = \left( \int_{-1}^1 \frac{1}{\ld - z} d E_{R;\ld} \vec e_k, \vec e_j \right)_{L^2(M)} =
   \left( ( R - zE_{L^2(M)} )^{-1} \vec e_k, \vec e_j \right)_{L^2(M)} = $$
$$ =  \left( \left( \frac{2}{b-a} Q - \frac{a+b}{b-a} E_{L^2(M)} - zE_{L^2(M)} \right)^{-1} \vec e_k,
   \vec e_j \right)_{L^2(M)} = $$
$$ = \int_a^b \left( \frac{2}{b-a} u - \frac{a+b}{b-a} - z \right)^{-1} d(E_u \vec e_k,\vec e_j)_{ L^2(M) } = $$
\begin{equation}
\label{f2_35}
= \int_{-1}^1 \frac{1}{\ld-z} d(E_{ \frac{b-a}{2}\ld + \frac{a+b}{2} } \vec e_k,\vec e_j)_{ L^2(M) },\
0\leq k,j\leq N-1,
\end{equation}
where $\{ \widehat E_\ld \}$ and $\{ \widehat E_{R;\ld} \}$ are left-continuous in $[-1,1)$, right-continuous
at the point $1$, constant outside $[-1,1]$, orthogonal
resolutions of unity of operators $\widehat B = URU^{-1}$ and $R$, respectively. Here $\{ E_\ld \}$ is
the orthogonal resolution of unity of $Q$, given by~(\ref{f2_7}).
By the Stieltjes-Perron inversion  formula we get
\begin{equation}
\label{f2_36}
( \widehat E_\ld x_k, x_j)_{\widehat H} =
(E_{ \frac{b-a}{2}\ld + \frac{a+b}{2} } \vec e_k, \vec e_j)_{L^2(M)},\qquad 0\leq k,j\leq N-1,
\end{equation}
for each $\ld\in [-1,1]$, such that $\ld$ is a point of continuity for $\widehat E_\ld$ and
$E_{ \frac{b-a}{2}\ld + \frac{a+b}{2} }$.
Using the change of variable we obtain that
\begin{equation}
\label{f2_37}
\widehat m_{k,j}(x) = ( P^{\widehat H}_H \widehat E_{ \frac{2}{b-a}x - \frac{a+b}{b-a} } x_k,x_j )_H =
(E_{x} \vec e_k, \vec e_j)_{L^2(M)},
\end{equation}
where $0\leq k,j\leq N-1$, for $x\in [a,b]$: $x$ is a point of continuity of
$\widehat E_{ \frac{2}{b-a}x - \frac{a+b}{b-a} }$ and $E_x$.
Using~(\ref{f2_7}) we can write
$$ (E_{x} \vec e_k, \vec e_j)_{L^2(M)} = m_{k,j}(x),\quad x\in (a,b), $$
and therefore
\begin{equation}
\label{f2_38}
\widehat m_{k,j} (x) = m_{k,j}(x),
\end{equation}
where $x\in [a,b]$: $x$ is a point of continuity of
$\widehat E_{ \frac{2}{b-a}x - \frac{a+b}{b-a} }$ and $E_x$.
Since matrix functions $\widehat M(x)$ and $M(x)$ are left-continuous in $(a,b)$, they
coincide in $(a,b)$. It remains to note that $\widehat M(a)=M(a)=0$, and $\widehat M(b)=M(b)=S_0$, to obtain
$$ \widehat M(x) = M(x),\qquad x\in [a,b]. $$
Consequently, all solutions of the truncated moment problem are generated by
left-continuous in $[-1,1)$, right-continuous at the point $1$, constant outside $[-1,1]$
sc-spectral functions of $B$.

\noindent
It remains to prove that different sc-spectral functions of the operator $B$ produce different
solutions of the moment problem~(\ref{f1_1}).
Suppose to the contrary that two different
left-continuous in $[-1,1)$, right-continuous at the point $1$, constant outside $[-1,1]$
sc-spectral functions produce the same solution of
the moment problem. That means that
there exist two self-adjoint contractions
$B_j\supseteq B$, in Hilbert spaces $H_j\supseteq H$, such that
\begin{equation}
\label{f2_39}
P_{H}^{H_1} E_{1,\ld} \not= P_{H}^{H_2} E_{2,\ld},
\end{equation}
\begin{equation}
\label{f2_40}
(P_{H}^{H_1} E_{1,\ld} x_k,x_j)_H = (P_{H}^{H_2} E_{2,\ld} x_k,x_j)_H,\qquad 0\leq k,j\leq N-1,\quad \ld\in[-1,1],
\end{equation}
where $\{ E_{n,\ld} \}_{\ld\in\R}$ are left-continuous in $[-1,1)$, right-continuous
at the point $1$, constant outside $[-1,1]$, orthogonal resolutions of unity of
operators $B_n$, $n=1,2$.
Set $L_N := \Lin\{ x_k \}_{k=0,N-1}$. By linearity we get
\begin{equation}
\label{f2_41}
(P_{H}^{H_1} E_{1,\ld} x,y)_H = (P_{H}^{H_2} E_{2,\ld} x,y)_H,\qquad x,y\in L_N,\quad \ld\in[-1,1].
\end{equation}
Denote by $R_{n,\ld}$ the resolvent of $B_n$, and set $\mathbf R_{n,\ld} := P_{H}^{H_n} R_{n,\ld}$, $n=1,2$.
From~(\ref{f2_41}),(\ref{f2_30}) it follows that
\begin{equation}
\label{f2_42}
(\mathbf R_{1,\ld} x,y)_H = (\mathbf R_{2,\ld} x,y)_H,\qquad x,y\in L_N,\quad \ld\in\CC\backslash\R.
\end{equation}
Choose an arbitrary $z\in\CC\backslash\R$ and consider the space $H_z :=
\overline{ (B-zE_H) H_a }$.
Since
$$ R_{j,z} (B-zE_H) x = (B_j - z E_{H_j} )^{-1} (B_j - z E_{H_j}) x = x,\qquad x\in H_a = D(B),$$
we get
\begin{equation}
\label{f2_43}
R_{1,z} u = R_{2,z} u \in H,\qquad u\in H_z;
\end{equation}
\begin{equation}
\label{f2_44}
\mathbf R_{1,z} u = \mathbf R_{2,z} u,\qquad u\in H_z,\ z\in\CC\backslash\R.
\end{equation}
We can write
$$ (\mathbf R_{n,z} x, u)_H = (R_{n,z} x, u)_{H_n} = ( x, R_{n,\overline{z}}u)_{H_n} =
( x, \mathbf R_{n,\overline{z}} u)_H,\ x\in H,\ u\in H_{\overline z}, $$
\begin{equation}
\label{f2_45}
n=1,2,
\end{equation}
and therefore we get
\begin{equation}
\label{f2_46}
(\mathbf R_{1,z} x,u)_H = (\mathbf R_{2,z} x,u)_H,\qquad x\in H,\ u\in H_{\overline z}.
\end{equation}
Choose an arbitrary $u\in H$, $u = \sum_{k=0}^{dN+N-1} c_k x_k$, $c_k\in\CC$.
Consider the following system of linear equations:
\begin{equation}
\label{f2_47}
-\left( \frac{a+b}{b-a} + z \right) d_k = c_k,\qquad  k=0,1,...,N-1;
\end{equation}
\begin{equation}
\label{f2_48}
\frac{2}{b-a} d_{k-N} - \left( \frac{a+b}{b-a} + z \right) d_k = c_k,\qquad  k=N,N+1,\ldots,dN+N-1;
\end{equation}
where $\{ d_k \}_{k=0}^{dN+N-1}$ are unknown complex numbers, $z\in\CC\backslash\R$ is a fixed parameter,
$a,b$ are from~(\ref{f1_1}).
Set
$$ d_k = 0,\qquad k=dN,dN+1,...,dN+N-1; $$
\begin{equation}
\label{f2_49}
d_{k-N} = \frac{b-a}{2} \left( \left( \frac{a+b}{b-a} + z \right) d_k + c_k \right),\qquad  k=dN+N-1,dN+N-2,...,N;
\end{equation}
For such numbers $\{ d_k \}_{k=0}^{dN+N-1}$, all equations in~(\ref{f2_48}) are satisfied.
Equations~(\ref{f2_47}) are not necessarily satisfied. Set
$v = \sum_{k=0}^{dN+N-1} d_k x_k = \sum_{k=0}^{dN-1} d_k x_k$. Notice that $v\in H_a = D(B)$.
We can write
$$ (B-zE_H) v = \left( \frac{2}{b-a} A - \frac{a+b}{b-a} E_H - zE_H \right) v = $$
$$ = \sum_{k=0}^{dN-1} d_k \left(
\frac{2}{b-a} x_{k+N} - \left( \frac{a+b}{b-a} + z \right) x_k
\right) = $$
$$ = \sum_{k=0}^{dN+N-1} \left( \frac{2}{b-a} d_{k-N} - \left( \frac{a+b}{b-a} + z \right) d_k \right) x_k, $$
where $d_{-1}=d_{-2}=...=d_{-N}=0$.
By the construction of $d_k$ we have
$$ (B-zE_H) v - u = \sum_{k=0}^{N-1} \left( - \left( \frac{a+b}{b-a} + z \right) d_k - c_k \right) x_k; $$
\begin{equation}
\label{f2_50}
u = (B-zE_H) v + \sum_{k=0}^{N-1} \left( \left( \frac{a+b}{b-a} + z \right) d_k + c_k \right) x_k,\qquad
u\in H,\ z\in\CC\backslash\R.
\end{equation}
By~(\ref{f2_50}) an arbitrary element $y\in H$ can be represented as $y=y_{ \overline{z} } + y'$,
$y_{ \overline{z} }\in H_{ \overline{z} }$, $y'\in L_N$.
Using~(\ref{f2_42}) and~(\ref{f2_46})  we get
$$ (\mathbf R_{1,z} x,y)_H = (\mathbf R_{1,z} x, y_{ \overline{z} } + y')_H =
(\mathbf R_{2,z} x, y_{ \overline{z} } + y')_H = (\mathbf R_{2,z} x,y)_H,\ x\in L_N,\ y\in H. $$
Thus, we obtain
\begin{equation}
\label{f2_51}
\mathbf R_{1,z} x = \mathbf R_{2,z} x,\qquad x\in L_N,\ z\in\CC\backslash\R.
\end{equation}
For an arbitrary $x\in L$, $x=x_z + x'$, $x_z\in H_z$, $x'\in L_N$, using
relations~(\ref{f2_44}),(\ref{f2_51}) we obtain
\begin{equation}
\label{f2_52}
\mathbf R_{1,z} x = \mathbf R_{1,z} (x_z + x') =
\mathbf R_{2,z} (x_z + x') = \mathbf R_{2,z} x,\qquad x\in L,\ z\in\CC\backslash\R,
\end{equation}
and
\begin{equation}
\label{f2_53}
\mathbf R_{1,z} x = \mathbf R_{2,z} x,\qquad x\in H,\ z\in\CC\backslash\R.
\end{equation}
By~(\ref{f2_30}) that means that the corresponding sc-spectral functions coincide and we obtain a
contradiction.
$\Box$

We shall recall some known facts about sc-resolvents (see~\cite{cit_9000_KO},\cite{cit_10000_KO}).
Let $A$ be a Hermitian contraction in a Hilbert space $H$ with a non-dense closed domain $\mathcal D=D(A)$,
and $\mathcal R = H \ominus \mathcal D$. A set of all self-adjoint extensions of $A$ in $H$, which are
contractions, we denote by $\mathcal{B}_H(A)$.
A set of all self-adjoint extensions of $A$ in a Hilbert space $\widetilde H\supseteq H$, which are
contractions, we denote by $\mathcal{B}_{\widetilde H} (A)$.
The set $\mathcal{B}_H(A)$ is non-empty. Moreover, there are a "minimal"$~$ element $A^\mu$ and
a "maximal"$~$ element $A^M$ in this set, such that $\mathcal{B}_H(A)$ coincides with the operator
segment
\begin{equation}
\label{f2_54}
A^\mu \leq \widetilde A\leq A^M.
\end{equation}
In the case $A^\mu = A^M$ the set $\mathcal{B}_H(A)$ consists of a unique element. This case is called
{\it determinate}.

\noindent
In the case $A^\mu \not= A^M$ the set $\mathcal{B}_H(A)$ consists of an infinite number of elements.
This case is called {\it indeterminate}.

\noindent
The case $A^\mu x \not= A^M x$, $x\in \mathcal{R}\backslash\{ 0 \}$, is called
{\it completely indeterminate}.
The indeterminate case can be always reduced to the completely indeterminate. If
$\mathcal{R}_0 = \{ x\in \mathcal{R}:\ A^\mu x = A^M x\}$, we can set
\begin{equation}
\label{f2_55}
A_e x = A x,\ x\in \mathcal{D};\quad A_e x = A^\mu x,\ x\in \mathcal{R}_0.
\end{equation}
The sets of generalized sc-resolvents for $A$ and for $A_e$ coincide (\cite[p. 1039]{cit_10000_KO}).

\noindent
Elements of $\mathcal{B}_H(A)$ are canonical (i.e. inside $H$) extensions of $A$ and their
resolvents are called canonical sc-resolvents of $A$.
On the other hand, elements  of $\mathcal{B}_{\widetilde H} (A)$ for $\widetilde H\supseteq H$
generate generalized sc-resolvents of $A$. The set of all generalized sc-resolvents we denote by
$\mathcal{R}^c(A)$.

\noindent
Set
\begin{equation}
\label{f2_56}
C = A^M - A^\mu,
\end{equation}
\begin{equation}
\label{f2_57}
Q_\mu (z) = \left.\left( C^{\frac{1}{2}} R^\mu_z C^{\frac{1}{2}} + E_H \right)\right|_{\mathcal{R}},\qquad
z\in\CC\backslash[-1,1],
\end{equation}
where $R^\mu_z = (A^\mu - zE_H)^{-1}$.

\noindent
An operator-valued function $k(z)$ with values in $[H]$ belongs to the class $R_{\mathcal{R}}[-1,1]$ if

\noindent
1) $k(z)$ is analytic in $z\in\CC\backslash[-1,1]$ and
$$ \frac{ \img k(z) }{ \img z }\leq 0,\qquad z\in\CC:\ \img z\not= 0; $$
2) For $z\in \R\backslash[-1,1]$, $k(z)$ is a self-adjoint positive contraction.
\begin{thm}
\label{t2_4} (\cite[p. 1053]{cit_10000_KO}).
The following equality:
\begin{equation}
\label{f2_58}
\widetilde R_z^c = R^\mu_z - R^\mu_z C^{\frac{1}{2}}  \left(
E_H + (Q_\mu(z)-E) k(z)
\right)^{-1}
C^{\frac{1}{2}} R^\mu_z,
\end{equation}
where $k(z)\in R_{\mathcal{R}}[-1,1]$, $\widetilde R_z^c\in \mathcal{R}^c(A)$,
establishes a one-to-one correspondence between the set $R_{\mathcal{R}}[-1,1]$ and the set $\mathcal{R}^c(A)$.

Moreover, the canonical resolvents correspond in~(\ref{f2_58}) to the constant functions
$k(z)\equiv K$, $K\in [0,E_{\mathcal{R}}]$.
\end{thm}
Comparing Theorem~\ref{t2_3} and Theorem~\ref{t2_4} we obtain the following result.
\begin{thm}
\label{t2_5}
Let a matrix moment problem~(\ref{f1_1}) with $l=2d$, $d\in\N$, be given.
Suppose that conditions~(\ref{f1_5}),(\ref{f2_5}) are true.
Let the operator $B$ be defined by~(\ref{f2_6}). The following statements are true:

\noindent
1) If $B^\mu = B^M$, then the moment problem~(\ref{f1_1}) has a unique solution. This solution
is given by
\begin{equation}
\label{f2_59}
M(x) = (m_{j,n}(x))_{j,n=0}^{N-1},\quad
m_{j,n}(x) = ( E^\mu_{ \frac{2}{b-a}x - \frac{a+b}{b-a} } x_j,x_n )_H,\qquad
0\leq j,n\leq N-1,
\end{equation}
where $\{ E^\mu_z \}$ is the left-continuous in $[-1,1)$, right-continuous
at the point $1$, constant outside $[-1,1]$, orthogonal resolution of unity of the operator $A^\mu$.

\noindent
2) If $B^\mu \not= B^M$, we define the extended operator $B_e$ (see the construction in~(\ref{f2_55})).
All solutions of the moment problem have the following form
\begin{equation}
\label{f2_60}
M(x) = (m_{j,n}(x))_{j,n=0}^{N-1},\quad
m_{j,n}(x) = ( \mathbf E_{ \frac{2}{b-a}x - \frac{a+b}{b-a} } x_j,x_n )_H,\qquad
0\leq j,n\leq N-1,
\end{equation}
where $\mathbf E_z$ is a left-continuous in $[-1,1)$, right-continuous
at the point $1$, constant outside $[-1,1]$ functions given by
$$ (R^\mu_z - R^\mu_z C^{\frac{1}{2}}  \left(
E_H + (Q_\mu(z)-E) k(z)
\right)^{-1}
C^{\frac{1}{2}} R^\mu_z f,g)_H =
\int_\R \frac{1}{\ld - z} d( \mathbf E_\ld f,g)_H, $$
\begin{equation}
\label{f2_61}
f,g\in H,\ z\in\CC\backslash\R,
\end{equation}
where $k(z)\in R_{\mathcal{R}}[-1,1]$. Here $\mathcal R = H \ominus D(B_e)$, and $C$,$Q_\mu(z)$,$R^\mu_z$
are defined by~(\ref{f2_56}),(\ref{f2_56}) with $A=B_e$.
Moreover, the correspondence between all solutions of the moment problem and $k(z)\in R_{\mathcal{R}}[-1,1]$
is one-to-one.
\end{thm}
P r o o f. It remains to consider the case 1). In this case all self-adjoint contractions
$\widetilde B\supseteq B$ in a Hilbert space $\widetilde H\supseteq H$ coincide on $H$ with $B^\mu$,
see~\cite[p. 1039]{cit_10000_KO}. Thus, the corresponding sc-spectral functions are
spectral functions of the self-adjoint operator $B^\mu$. However, a self-adjoint operator
has a unique (normalized) spectral function. Thus, a set of sc-spectral functions of $B$ consists of
a unique element. This element is a spectral function of $B^\mu$.
$\Box$

Consider a matrix moment problem~(\ref{f1_1}) with $l=0$. In this case the necessary and sufficient
condition of solvability is
\begin{equation}
\label{f2_62}
S_0 \geq 0.
\end{equation}
The necessity is obvious. On the other hand, if relation~(\ref{f2_62}) is true, we can choose
\begin{equation}
\label{f2_63}
M(x) = \frac{x-a}{b-a} S_0,\qquad x\in[a,b].
\end{equation}
This function is a solution of the moment problem.
A set of all solutions consists of non-decreasing matrix functions $M(x)$ on $[a,b]$, left-continuous in
$(a,b)$, with the boundary conditions $M(a)=0$, $M(b)=S_0$.
\section{The case of an even number of given moments: solvability and a description  of solutions.}
Consider the matrix moment problem~(\ref{f1_1}) with $l=2d+1$, $d\in\Z_+$.
\begin{thm}
\label{t3_1}
Let a matrix moment problem~(\ref{f1_1}) with $l=2d+1$, $d\in\Z_+$, be given.
The moment problem has a solution if and only if
\begin{equation}
\label{f3_1}
\Gamma_{d}\geq 0;\ \widetilde\Gamma_{d}\geq 0,
\end{equation}
and there exist matrix solutions $X,Y$ of matrix equations
\begin{equation}
\label{f3_2}
\Gamma_{d}X = \left(\begin{array}{cccc} S_{d+1}\\
S_{d+2}\\
\vdots\\
S_{2d+1}\end{array}\right),\quad
\widetilde \Gamma_{d}Y = \left(\begin{array}{cccc} -ab S_{d}+(a+b)S_{d+1}-S_{d+2}\\
-ab S_{d+1}+(a+b)S_{d+2}-S_{d+3}\\
\vdots\\
-ab S_{2d-1}+(a+b)S_{2d}-S_{2d+1}\end{array}\right),
\end{equation}
and for these solutions $X$, $Y$ the following relation is true:
\begin{equation}
\label{f3_3}
X^*\Gamma_{d}X \leq -abS_{2d} + (a+b)S_{2d+1} - Y^* \widetilde\Gamma_{d}Y.
\end{equation}
\end{thm}
P r o o f.
Consider a matrix moment problem~(\ref{f1_1}) with $l=2d+1$, $d\in\Z_+$.
It has a solution if and only if the moment problem with an odd number of moments
\begin{equation}
\label{f3_4}
\int_a^b x^n dM(x) = S_n,\ n=0,1,...,2d+1;\quad \int_a^b x^{2d+2} dM(x) = S_{2d+2},
\end{equation}
with some complex $(N\times N)$ matrix $S_{2d+2}$ has a solution.
By~(\ref{f1_5}) the solvability of the moment problem~(\ref{f3_4}) is equivalent to the matrix inequalities
\begin{equation}
\label{f3_5}
\Gamma_{d+1}\geq 0,\ \widetilde\Gamma_{d+1}\geq 0.
\end{equation}
If we apply to the latter inequalities the well-known block-matrix lemma (e.g.~\cite[p. 223]{cit_15000_Z}),
we obtain that solvability of~(\ref{f3_4}) is equivalent to the condition~(\ref{f3_1}), existence
of solutions $X,Y$ of~(\ref{f3_2}) and inequalities
\begin{equation}
\label{f3_6}
S_{2d+2}\geq X^*\Gamma_d X,\quad S_{2d+2}\leq -abS_{2d}+(a+b)S_{2d+1} - Y^*\widetilde\Gamma_d Y.
\end{equation}
Consequently, we get that the statement of the Theorem is true.
$\Box$

\begin{thm}
\label{t3_2}
Let a matrix moment problem~(\ref{f1_1}) with $l=2d+1$, $d\in\Z_+$, be given. Suppose that
conditions~(\ref{f3_1}),(\ref{f3_2}) and~(\ref{f3_3}) are true.
For an arbitrary $S_{2d+2}\in [X^*\Gamma_{d}X \leq -abS_{2d} + (a+b)S_{2d+1} - Y^* \widetilde\Gamma_{d}Y]$,
the set of all solutions of the matrix moment  problem with an odd number of moments~(\ref{f3_4}) we
denote by $\mathcal{V}(S_{2d+2})$.
The set $\mathcal{V}$ of all solutions of the moment problem~(\ref{f1_1}) is given by the formula
\begin{equation}
\label{f3_7}
\mathcal{V} = \cup_{S\in [X^*\Gamma_{d}X , -abS_{2d} + (a+b)S_{2d+1} - Y^* \widetilde\Gamma_{d}Y]} \mathcal{V}(S).
\end{equation}
The sets $\mathcal{V}(S)$ in~(\ref{f3_7}) for different $S$ do not intersect. Each set $\mathcal{V}(S)$
in~(\ref{f3_7}) is parameterized by virtue of Theorem~\ref{t2_5}.
\end{thm}
P r o o f.
The proof of~(\ref{f3_7}) follows obviously from the above considerations in the proof of Theorem~\ref{t3_1}.
The sets $\mathcal{V}(S)$ in~(\ref{f3_7}) for different $S$ do not intersect since
the solutions in different sets have different $(2d+2)$-th moments.
$\Box$

\renewcommand{\refname}{\Large \rm\bf

\begin{center}
\bf Truncated matricial moment problems on a finite interval: the operator approach.
\end{center}
\begin{center}
\bf
S.M. Zagorodnyuk
\end{center}
In this paper we obtain a description of all solutions of truncated matricial moment problems on a finite interval
in a general case (no conditions besides solvability are assumed). We use the basic results of M.G.~Krein and
I.E.~Ovcharenko about generalized sc-resolvents of Hermitian contractions.

MSC 2000: 44A60, 30E05.

Key words: moment problem, generalized resolvent, contraction.


\centerline{References.}}
\vskip1cm

\begin{thebibliography}{99}
\bibitem{cit_1000_KN} {\it M.G. Krein, A.A. Nudelman}. The Markov moment problem and extremal
problems. Nauka. Moscow (1973). 552p. (Russian).
\bibitem{cit_2000_CDFK} {\it A.E.~Choque Rivero, Yu.M.~Dyukarev, B.~Fritzsche and B.~Kirstein}.
A truncated matricial moment problem on a finite interval. In: Interpolation, Schur functions and
moment problems (Eds.: D.~Alpay and I.~Gohberg), Operator Theory: Advances and Applications,
Birkh\"auser, Basel-Boston-Berlin. (2006), v.165, p.121-173.
\bibitem{cit_3000_CDFK} {\it A.E.~Choque Rivero, Yu.M.~Dyukarev, B.~Fritzsche and B.~Kirstein}.
A truncated matricial moment problem on a finite interval. The case of an odd number of
prescribed moments. In: System theory, the Schur algorithm and multidimensional analysis
(Eds.: D.~Alpay, V.~Vinnikov and I.~Gohberg), Operator Theory: Advances and Applications,
Birkh\"auser, Basel-Boston-Berlin. (2007), v.176, p.99-164.
\bibitem{cit_4000_N} {\it M.A. Neumark}. Spectral functions of a symmetric
operator. Izvestiya AN SSSR (1940), v.4, p. 277-318. (Russian).
\bibitem{cit_5000_N} {\it M.A. Neumark}. On spectral functions of a symmetric
operator. Izvestiya AN SSSR (1943), v.7, p. 285-296. (Russian).
\bibitem{cit_6000_KK} {\it M.G. Krein, M.A. Krasnoselskiy}.
Basic theorems on an extension of Hermitian operators and some their applications to
the theory of orthogonal polynomials and the moment  problem.
Uspehi matem. nauk (1947), v.3(19), p.60-106. (Russian).
\bibitem{cit_7000_A}
{\it N. I. Akhiezer}. Classical moment problem.
Fizmatlit. Moscow. (1961). 312p. (Russian).
\bibitem{cit_8000_Z} {\it S.M. Zagorodnyuk}. Positive definite kernels satisfying difference equations.
Methods of funct. analysis and topology (2010), v.1.
\bibitem{cit_9000_KO} {\it M.G. Krein, I.E. Ovcharenko}.
On the theory of generalized resolvents of
nondensely defined Hermitian contractions.
Dopovidi AN URSR (1976), seriya A, No.10, p.881-884. (Ukrainian).
\bibitem{cit_10000_KO} {\it M.G. Krein, I.E. Ovcharenko}.
On Q-functions and sc-resolvents of
nondensely defined Hermitian contractions.
Sibirskiy matem. zhurnal (1977), v.XVIII, No.5, p.1032-1056. (Russian).
\bibitem{cit_11000_K} {\it M.G. Krein}.
The theory of self-adjoint extensions of semi-bounded Hermitian transformations and its
applications. I.
Matem. sbornik (1947), v.20(62), No.3, p.431-495. (Russian).
\bibitem{cit_11500_MM} {\it M. M. Malamud, S. M. Malamud}.
Operator measures in a  Hilbert space.
Algebra i analiz (2003), v.15, No.3, p.~1-52.
\bibitem{cit_12000} {\it V.S. Korolyuk (ed.)}.
Handbook on the theory of  probability and mathematical statistics.
Naukova Dumka. Kiev (1978). 584p. (Russian).
\bibitem{cit_13000}
{\it Ju. M. Berezanskii}.
Expansions in eigenfunctions of selfadjoint operators.
Amer. Math. Soc., Providence, RI (1968). (Russian edition: Naukova Dumka, Kiev, 1965).
\bibitem{cit_13500_S}
{\it A.V. Shtraus}.
Generalized resolvents of symmetric operators.
Izvestiya AN SSSR (1954), v.18, p.51-86. (Russian).
\bibitem{cit_14000_AG}
{\it N. I. Akhiezer, I. M. Glazman}.
Theory of linear operators in a Hilbert space.
Gos. izdat. teh.-teor. lit. Moscow, Leningrad (1950). 484p. (Russian).
\bibitem{cit_15000_Z}
{\it V. A. Zolotarev}.
Analytic methods of spectral representations of non-selfadjoint and non-unitary operators.
KhNU. Kharkov (2003). 344p. (Russian).
\end{thebibliography}
\end{document}